\magnification=1200

          %%%%%%%%%%%%%%%%%%%%%%%%%%%%%%%%%%%%%%%%%%%%%%%%%%%%%%%%%%%
          %       FOP.tex FORTEX-compatible version of papmacs      %
          %   VERSION OF May 6, 1989; does not use the AMS-fonts    %
          %%%%%%%%%%%%%%%%%%%%%%%%%%%%%%%%%%%%%%%%%%%%%%%%%%%%%%%%%%%

\def\item{\vskip1.3pt\hang\textindent}% THIS REPLACES KNUTH'S DEF'N

% THIS
                                              %REPLACES KNUTH'S DEF'N

\tolerance=300 \pretolerance=200 \hfuzz=1pt \vfuzz=1pt

%\magnification=\magstep1
\hoffset 0cm            %all offsets went wrong
%\voffset=0.8true cm
\hsize=5.8 true in \vsize=9.5 true in

\def\rightheadline{\hfil\smc\lastname\hfil\tenbf\folio}
\def\leftheadline{\tenbf\folio\hfil\smc\lastname\hfil}
\headline={\ifodd\pageno\rightheadline\else\leftheadline\fi}
\newdimen\dimenone
\def\checkleftspace#1#2#3#4#5{%DIESER MACRO STAMMT VON APPELT
 \dimenone=\pagetotal
 \advance\dimenone by -\pageshrink   %testen ob Titel noch mit Gewalt auf Seite
                                                                          %geht
 \ifdim\dimenone>\pagegoal          %nacha tua nix-- gewoehnliche Outputroutine
   \else\dimenone=\pagetotal
        \advance\dimenone by \pagestretch
        \ifdim\dimenone<\pagegoal
          \dimenone=\pagetotal
          \advance\dimenone by#1         %addieren Skip vor Ueberschrift (=#1)
          \setbox0=\vbox{#2\parskip=0pt                %#2 ist gewaehlter Font
                       \hyphenpenalty=10000
                       \rightskip=0pt plus 5em
                       \noindent#3 \vskip#4}    %#3=Ueberschrift,#4=skip nachher
        \advance\dimenone by\ht0
        \advance\dimenone by 3\baselineskip
        \ifdim\dimenone>\pagegoal\vfill\eject\fi
          \else\eject\fi\fi}

\parindent=35pt
\mathsurround=1pt
\parskip=1pt plus .25pt minus .25pt
\normallineskiplimit=.99pt

\mathchardef\emptyset="001F % THIS REPLACES KNUTH'S DEFINITION

% USED FOR REAL PART OF COMPLEX NUMBERS
% USED FOR IMAGINARY PART OF COMPLEX NUMBERS
 % USED FOR IDENTITY FUNCTION
\def\Int{\mathop{\rm int}\nolimits}
% USED FOR TRACE OF MATRIX
%

%\def\u{\mathop{\rm u}\nolimits}

%\def\O{\mathop{\rm O}\nolimits}

\def\1{{\bf1}}\def\0{{\bf0}}

\def\({\bigl(}  \def\){\bigr)}
\def\<{\mathopen{\langle}}\def\>{\mathclose{\rangle}}

\def\Z{{\mathchoice{{\hbox{$\rm Z\hskip 0.26em\llap{\rm Z}$}}}%
{{\hbox{$\rm Z\hskip 0.26em\llap{\rm Z}$}}}%
{{\hbox{$\scriptstyle\rm Z\hskip 0.31em\llap{$\scriptstyle\rm Z$}$}}}{{%
\hbox{$\scriptscriptstyle\rm
Z$\hskip0.18em\llap{$\scriptscriptstyle\rm Z$}}}}}}

\def\F{{\mathchoice{\hbox{$\rm I\hskip-0.14em F$}}%
{\hbox{$\rm I\hskip-0.14em F$}}%
{\hbox{$\scriptstyle\rm I\hskip-0.14em F$}}%
{\hbox{$\scriptscriptstyle\rm I\hskip-0.10em F$}}}}

\def\R{{\mathchoice{\hbox{$\rm I\hskip-0.14em R$}}%
{\hbox{$\rm I\hskip-0.14em R$}}%
{\hbox{$\scriptstyle\rm I\hskip-0.14em R$}}%
{\hbox{$\scriptscriptstyle\rm I\hskip-0.10em R$}}}}

\def\Q{\QQ\,\,}

\def\.{{\cdot}}
\def\|{\Vert}
\def\ssk{\smallskip}
\def\msk{\medskip}
\def\bsk{\bigskip}
\def\giantskip{\vskip2\bigskipamount}

\def\giantbreak{\par \ifdim\lastskip<2\bigskipamount \removelastskip
         \penalty-400 \giantskip\fi}

\def\nin{\noindent}
\def\cen{\centerline}
\def\pagebreak{\vskip 0pt plus 0.0001fil\break}
\def\linebreak{\break}

\def\epsilon{\varepsilon}

\font\ninerm=cmr9 \font\eightrm=cmr8 \font\sixrm=cmr6
 \font\eightbf=cmbx8 \font\sixbf=cmbx6
 \font\eighti=cmmi8 \font\sixi=cmmi6
\font\ninesy=cmsy9 \font\eightsy=cmsy8 \font\sixsy=cmsy6
 \font\eightit=cmti8 
% SANS SERIF 10 POINT
 %SANS SERIF 10 POINT ITALIC
 \font\eightsl=cmsl8 
\font\eighttt=cmtt8
 %SLANTED TYPEWRITER 10 POINT
 %BOLD FACE MATH SYMBOLS 10 POINT
 %DUNHILL STYLE 10 POINT
 %SAN SERIF BOLD EXTENDED 10 POINT
 %USED FOR TITLES
 %USED FOR TITLES
\font\bfone=cmbx10 scaled\magstep1 %BOLDFACE AT MAGSTEP 1
 %BOLDFACE AT MAGSTEP 2
 %BOLDFACE AT MAGSTEP 3
\font\smc=cmcsc10 
 
scaled\magstep1 \font\small=cmcsc8

\def\no #1. {\bigbreak\vskip-\parskip\noindent\bf #1. \quad\rm}

\def\Proposition #1. {\checkleftspace{0pt}{\bf}{Theorem}{0pt}{}
\bigbreak\vskip-\parskip\noindent{\bf Proposition #1.} \quad\it}

\def\Theorem #1. {\checkleftspace{0pt}{\bf}{Theorem}{0pt}{}
\bigbreak\vskip-\parskip\noindent{\bf  Theorem #1.} \quad\it}
\def\Corollary #1. {\checkleftspace{0pt}{\bf}{Theorem}{0pt}{}
\bigbreak\vskip-\parskip\nin{\bf Corollary #1.} \quad\it}
\def\Lemma #1. {\checkleftspace{0pt}{\bf}{Theorem}{0pt}{}
\bigbreak\vskip-\parskip\noindent{\bf  Lemma #1.}\quad\it}

\def\Definition #1. {\checkleftspace{0pt}{\bf}{Theorem}{0pt}{}
\rm\bigbreak\vskip-\parskip\noindent{\bf Definition #1.} \quad}

\def\Remark #1. {\checkleftspace{0pt}{\bf}{Theorem}{0pt}{}
\rm\bigbreak\vskip-\parskip\noindent{\bf Remark #1.}\quad}

\def\Exercise #1. {\checkleftspace{0pt}{\bf}{Theorem}{0pt}{}
\rm\bigbreak\vskip-\parskip\noindent{\bf Exercise #1.} \quad}

\def\Example #1. {\checkleftspace{0pt}{\bf}{Theorem}{0pt}{}
\rm\bigbreak\vskip-\parskip\noindent{\bf Example #1.}\quad}
\def\Examples #1. {\checkleftspace{0pt}{\bf}{Theorem}{0pt}
\rm\bigbreak\vskip-\parskip\noindent{\bf Examples #1.}\quad}

\newcount\problemnumb \problemnumb=0
\def\Problem{\global\advance\problemnumb by 1\bigbreak\vskip-\parskip\noindent
{\bf Problem \the\problemnumb.}\quad\rm }

\def\Proof#1.{\rm\par\ifdim\lastskip<\bigskipamount\removelastskip\fi\smallskip
            \noindent {\bf Proof.}\quad}

\nopagenumbers

\def\author{}
\def\lastname{}
\def\thanks#1{\footnote*{\eightrm#1}}
\def\title{}

\def\nonumbers{\def\leftheadline{\hfil} \def\rightheadline{\hfil}}

\def\lastname{}
\def\h{{\textstyle{1\over2}}}

\def\he{{1\over2}}

\def\ep{\epsilon}

\def\text{\textstyle}
\def\disp{\displaystyle}
\def\d{{\,\rm d}}

\def\and{{\rm and }}

\expandafter\edef\csname amssym.def\endcsname{%
       \catcode`\noexpand\@=\the\catcode`\@\space}
%  Set the catcode to 11 for use in private control sequence names.
\catcode`\@=11
\def\undefine#1{\let#1\undefined}
\def\newsymbol#1#2#3#4#5{\let\next@\relax
 \ifnum#2=\@ne\let\next@\msafam@\else
 \ifnum#2=\tw@\let\next@\msbfam@\fi\fi
 \mathchardef#1="#3\next@#4#5}
\def\mathhexbox@#1#2#3{\relax
 \ifmmode\mathpalette{}{\m@th\mathchar"#1#2#3}%
 \else\leavevmode\hbox{$\m@th\mathchar"#1#2#3$}\fi}
\def\hexnumber@#1{\ifcase#1 0\or 1\or 2\or 3\or 4\or 5\or 6\or 7\or 8\or
 9\or A\or B\or C\or D\or E\or F\fi}

%Loading fontfiles `eufm' and `msbm'

\font\tenmsb=msbm10 \font\sevenmsb=msbm7 \font\fivemsb=msbm5
\newfam\msbfam
\textfont\msbfam=\tenmsb\scriptfont\msbfam=\sevenmsb
\scriptscriptfont\msbfam=\fivemsb \edef\msbfam@{\hexnumber@\msbfam}
\def\Bbb#1{{\fam\msbfam\relax#1}}

\font\teneufm=eufm10 \font\seveneufm=eufm7 \font\fiveeufm=eufm5
\newfam\eufmfam
\textfont\eufmfam=\teneufm \scriptfont\eufmfam=\seveneufm
\scriptscriptfont\eufmfam=\fiveeufm

\catcode`@=11 %Set the catcode to 11 for use in private control sequence names.

\expandafter\edef\csname amssym.def\endcsname{%
       \catcode`\noexpand\@=\the\catcode`\@\space}
\font\eightmsb=msbm8 \font\sixmsb=msbm6 \font\fivemsb=msbm5
\font\eighteufm=eufm8 \font\sixeufm=eufm6 \font\fiveeufm=eufm5
\newskip\ttglue
\def\eightpoint{\def\rm{\fam0\eightrm}%
  \textfont0=\eightrm \scriptfont0=\sixrm \scriptscriptfont0=\fiverm
  \textfont1=\eighti \scriptfont1=\sixi \scriptscriptfont1=\fivei
  \textfont2=\eightsy \scriptfont2=\sixsy \scriptscriptfont2=\fivesy
  \textfont3=\tenex \scriptfont3=\tenex \scriptscriptfont3=\tenex
\textfont\eufmfam=\eighteufm \scriptfont\eufmfam=\sixeufm
\scriptscriptfont\eufmfam=\fiveeufm \textfont\msbfam=\eightmsb
\scriptfont\msbfam=\sixmsb \scriptscriptfont\msbfam=\fivemsb
  \def\it{\fam\itfam\eightit}%
  \textfont\itfam=\eightit
  \def\sl{\fam\slfam\eightsl}%
  \textfont\slfam=\eightsl
  \def\bf{\fam\bffam\eightbf}%
  \textfont\bffam=\eightbf \scriptfont\bffam=\sixbf
   \scriptscriptfont\bffam=\fivebf
  \def\tt{\fam\ttfam\eighttt}%
  \textfont\ttfam=\eighttt
  \tt \ttglue=.5em plus.25em minus.15em
  \normalbaselineskip=9pt
  \def\MF{{\manual opqr}\-{\manual stuq}}%
  \let\big=\eightbig
  \setbox\strutbox=\hbox{\vrule height7pt depth2pt width\z@}%
  \normalbaselines\rm}
\def\eightbig#1{{\hbox{$\textfont0=\ninerm\textfont2=\ninesy
  \left#1\vbox to6.5pt{}\right.\n@space$}}}

%  Restore the catcode value for @ that was previously saved.

\csname amssym.def\endcsname

% end of \symb1.tex

\def\la{\lambda}
\def\al{\alpha}
\def\be{\beta}

\def\om{\omega}
\def\({\left(}
\def\){\right)}
\def\for{\qquad \hbox{for}\ }
\def\eq{\eqalign}

\def\O#1{O\(#1\)}
\def\abs#1{\left| #1 \right|}

\def\norm#1{\left\Vert #1 \right\Vert}

\def\klein{\eightpoint \def\smc{\small} \baselineskip=9pt}

\def\fn#1#2{{\parindent=0.7true cm
\footnote{$^{(#1)}$}{{\klein  #2}}}}

\font\boldmas=msbm10                  %%
\def\Bbb#1{\hbox{\boldmas #1}}        %%
\def\Z{{\Bbb Z}}                        %%
\def\Q{{\Bbb Q}}
\def\R{{\Bbb R}}
\def\F{{\Bbb F}}

                  %%
        %%

%%%%%%%%%%%%%%%%%%%%%%%%%%%%%%%%%%%%%%%%%%%%%%%%%%%%%%%%%%%%%%%%%%%%%%
\font\eightrm=cmr8 \long\def\fussnote#1#2{{\baselineskip=9pt
\setbox\strutbox=\hbox{\vrule height 7pt depth 2pt width 0pt}%
\eightrm \footnote{#1}{#2}}}
%%%%%%%%%%%%%%%%%%%%%%%%%%%%%%%%%%%%%%%%
%% zum Verkleinern (Summationsgrenzen, Folgenindex, e.a.)%%
%%%%%%%%%%%%%%%%%%%%%%%%%%%%%%%%%%%%%%%%
\font\boldmasi=msbm10 scaled 700      %%
\def\Bbbi#1{\hbox{\boldmasi #1}}      %%
\font\boldmas=msbm10                  %%
\def\Bbb#1{\hbox{\boldmas #1}}        %%
\def\Zi{{\Bbbi Z}}                      %%
\def\Pi{{\Bbbi P}}                      %%
\def\Ri{{\Bbbi R}}

                      %%
%%%%%%%%%%%%%%%%%%%%%%%%%%%%%%%%%%%%%%%%

                        %%

\def\dint #1 {% Doppelintegral; #1 Integrationsbereich
\quad  \setbox0=\hbox{$\disp\int\!\!\!\int$}
  \setbox1=\hbox{$\!\!\!_{#1}$}
  \vtop{\hsize=\wd1\centerline{\copy0}\copy1} \quad}

\def\drint #1 {% Dreifachintegral; #1 Integrationsbereich
\qquad  \setbox0=\hbox{$\disp\int\!\!\!\int\!\!\!\int$}
  \setbox1=\hbox{$\!\!\!_{#1}$}
  \vtop{\hsize=\wd1\centerline{\copy0}\copy1}\qquad}

\def\frac#1#2{{#1\over #2}}

\def\date{\the\day.~\the\month.~\the\year}

\def\mod{\,{\rm mod}\,}
\def\klein{\eightpoint \def\smc{\small} }

\def\frac#1#2{{#1\over#2}}
\def\Int{\int\limits}

\def\vol{{\rm vol}}

\nonumbers

\hsize=16.4true cm     \vsize=23.3true cm

\parindent=0cm

\def\eqno{\leqno}
\def\l{\ell}
\def\DS{\sum_{1\le h\le U}\ \sum_{1\le m\le u}}
\def\Sh{\sum_{h=1}^{[U]}}
\def\alh{\al_{h,[U]}}  \def\beh{\be_{h,[U]}} \def\gah{\gamma_{h,[U]}}
\def\Mj{{\cal M}_j}
\def\J{{\cal J}}
\def\c{{\cal C}}
\def\D{{\cal D}}
\def\E{{\cal E}}
\def\Fe{{\cal F}}
\def\odd{\ {\rm odd}}
\def\b#1{{\bf #1}}
\def\rr{r_1\cdot\dots\cdot r_\l}

\vbox{\vskip 1.5true cm}

\footnote{}{\klein{\it Mathematics Subject Classification }
(2000): 11N37, 35P20, 58J50, 11P21.\par }

\cen{{\bfone A lower bound for the error term in Weyl's law}} \msk
\cen{{\bfone for certain Heisenberg manifolds, II}}\bsk \cen{{\bf
 Werner Georg Nowak }\fn{*}{The author gratefully
acknowledges support from the Austrian Science Fund (FWF) under
project Nr.~P20847-N18.} {\bf(Vienna)} }

\vbox{\vskip 1.2true cm}

{\klein{\bf Abstract. } This article is concerned with estimations
from below for the remainder term in Weyl's law for the spectral
counting function of certain rational $(2\l+1)$-dimensional
Heisenberg manifolds. Concentrating on the case of odd $\l$, it
continues the work done in part I [20] which dealt with even $\l$.}

\vbox{\vskip 1true cm}

{\bf 1.~Introduction. Weyl's law and Heisenberg manifolds. } Let
$M$ be a closed $n$-dimensional Riemannian manifold with a metric
$g$ and Laplace-Beltrami operator $\Delta$. Denote by $N(t)$ the
spectral counting function
$$ N(t):=\sum_{\lambda{\ \rm eigenvalue\ of\ }\Delta\atop
\lambda\le t} d(\lambda) $$ where $d(\lambda)$ is the dimension of
the eigenspace corresponding to $\lambda$, and $t$ is a large real
variable. According to a deep general theorem of L.~H\"{o}rmander [10],
$$ N(t) = {\vol(M)\over(4\pi)^{n/2}\,\Gamma(\h n+1)}\,t^{n/2} +
\O{t^{(n-1)/2}}\,, \eqno(1.1) $$ where the error term - in this
general setting - is best possible. Asymptotics like (1.1) and its
refinements for special manifolds are usually subsumed under the
notion of Weyl's law. \msk

Recently, the spectral theory of so-called Heisenberg manifolds has
attracted a lot of attention, possibly on the grounds of motivation
from quantum physics and the abstract theory of PDE's. To recall
basics, let $\l\ge1$ be an integer, and put\fn{1}{Bold face letters
will denote throughout elements of some space $\Ri^n$, resp., of
$\Zi^n$. They may be viewed also as $(1\times n)$-matrices ("row
vectors") where applicable.}
$$ \gamma(\b x,\b y,z) = \pmatrix{1 &\b x&z\cr ^t\b o_\l&I_\l& ^t\b
y\cr 0&\b o_\l&1\cr}\,,
$$ where $\b x, \b y \in \R^\l$, $z\in\R$,
$\b o_\l =(0,\dots,0)\in\R^\l$, $I_\l$ is the $(\l\times\l)$-unit
matrix, and $^t\cdot$ denotes transposition. Then the
$(2\l+1)$-dimensional Heisenberg group $H_\l$ is defined by $$ H_\l
= \{\gamma(\b x,\b y,z):\ \b x, \b y \in \R^\l,\ z\in\R\ \}\,,
\eqno(1.2) $$ with the usual matrix product. Further, for any
$\l$-tuple $\b r=(r_1,\dots,r_\l)\in\Z_+^\l$ with the property that
$r_j\mid r_{j+1}$ for all $j=1,\dots,\l-1$, let $\b r * \Z^\l:=
r_1\Z\times\dots\times r_\l\Z $, and define $$ \Gamma_{\b r} =
\{\gamma(\b x,\b y,z):\ \b x \in \b r*\Z^\l, \b y \in\Z^\l , z\in\Z\
\}\,. \eqno(1.3)$$ $\Gamma_{\b r}$ is a uniform discrete subgroup of
$H_\l$, i.e., the {\it Heisenberg manifold } $H_\l/\Gamma_{\b r}$ is
compact. According to Gordon and Wilson [6], Theorem 2.4, the
subgroups $\Gamma_{\b r}$ classify {\it all } uniform discrete
subgroups of $H_\l$ up to automorphisms: For every uniform discrete
subgroup $\Gamma$ of $H_\l$ there exists a unique $\l$-tuple $\b r$
and an automorphism of $H_\l$ which maps $\Gamma$ to $\Gamma_{\b
r}$. \ssk However, to get a "rational" or "arithmetic" Heisenberg
manifold - this latter expression being due to Petridis \& Toth [21]
- one has to make a quite particular choice of the
metric\fn{2}{Compare the discussion below concerning the bound (3.3)
which applies to "almost all" metrics $g$.} $g$. Following Petridis
\& Toth [21], Theorem 1.1, and also Zhai [24], we choose
$$ g_\l = \pmatrix{I_{2\l}& ^t\b o_{2\l}\cr \b o_{2\l}&2\pi}\,. \eqno(1.4)$$
The spectrum of the Laplace-Beltrami operator on $(H_\l/\Gamma_{\b
r},g_\l)$ has been analyzed in Gordon and Wilson [6], p.~259, and
also in Khosravi and Petridis [16], p.~3564. It consists of two
different classes ${\cal S}_I$ and ${\cal S}_{II}$, where ${\cal
S}_{I}$ is the spectrum of the Laplacian on the $2\l$-dimensional
torus $\R^{2\l}/\Z^{2\l}$, and
$$ {\cal S}_{II} = \{2\pi\(n_0^2+ n_0(2n_1+\l)\)\,:\
n_0\in\Z^+,\ n_1\in\Z_0^+ \ \}\,, \eqno(1.5)$$ with multiplicities (
= dimensions of corresponding eigenspaces)
$2n_0^\l\,\rr\,{n_1+\l-1\choose\l-1}$.  \bsk\msk

{\bf 2.~Lattice points in a circle. } The quantity $$
\sum_{\la\in{\cal S}_{II}\atop\la\le t} d(\la)\,,
$$ yields the major contribution to $N(t)$ for this rational
Heisenberg manifolds. Its asymptotic evaluation amounts to the
enumeration of the integer points $(n_0,n_1)\in(\Z_0^+, \Z^+)$ in
the planar domain $u^2+u(2v+\l)\le t/(2\pi)$, with the weights
$2n_0^\l\,\rr\,{n_1+\l-1\choose\l-1}$ as indicated. This observation
may be considered as {\it one } motivation to make reference to the
state-of-art with the Gaussian circle problem, the "ancestor and
prototype" of all planar lattice point problems. As a second link,
one may notice that $M=\R^2/\Z^2$, the 2-dimensional torus, is the
simplest example of a Riemannian manifold with a non-trivial
spectral theory: In fact\fn{3}{See also the detailed discussion in
part I of this work [20].}, the eigenvalues of the Laplacian on
$\R^2/\Z^2$ are given by $4\pi^2 k$, where $k$ ranges over all
nonnegative integers with $r(k)>0$, $r(k)$ denoting as usual the
number of ways to write $k$ as a sum of two squares of integers. The
corresponding multiplicities are given by $r(k)$, hence the spectral
counting function $N(t)$ now equals the number of lattice points in
an origin-centered compact circular disc of radius
$\sqrt{t}/(2\pi)$. \ssk For enlightening accounts on the history of
the Gaussian circle problem in textbook style, the reader may
consult the monographs of Kr\"{a}tzel [18], [19], and Huxley [11], along
with the recent quite comprehensive survey article [13]. The
sharpest upper bound for the {\it lattice point discrepancy } $P(x)$
of the compact unit circular disc $\D_0$, linearly dilated by a
large real parameter $x$, is nowadays due to Huxley [12] and reads
$$ P(x):= \#(x\D_0\cap\Z^2)-\pi x^2=\O{x^{131/208}(\log x)^{18637/8320}}\,.
\eqno(2.1)$$ It is usually conjectured that $ P(x)=\O{x^{1/2+\ep}}$
for every $\ep>0$. This is supported by Cram\'er's [3] classic
mean-square asymptotics $$ \Int_0^X (P(x))^2\d x \sim C\,X^2
\eqno(2.2)  $$ with an explicit constant $C>0$. Thus, roughly
speaking, $P(x)\ll x^{1/2}$ {\it in square-mean}, but it has been
known for a long time that there exist unbounded sequences of
$x$-values for which $P(x)$ attains "exceptionally large" values,
even of either sign: By 1961, the state-of-art in this direction was
that\fn{4}{Recall the usual $\Omega$-notation: For real functions
$F$ and $G>0$, and $*$ denoting either $+$ or $-$,
$F(x)=\Omega_*(G(x))$ means that $\lim\sup(*F(x)/G(x))>0$, as
$x\to\infty$. Further, $F(x)=\Omega(G(x))$ means that
$\lim\sup\abs{F(x)/G(x)}>0$.}
$$ P(x)=\Omega_-\(x^{1/2}(\log x)^{1/4}\)\eqno(2.3) $$ and $$
P(x)=\Omega_+\(x^{1/2}(\log_2 x\,\log_3 x)^{1/4}\)\,, \eqno(2.4)
$$ due to Hardy [9], resp., Gangadharan [5]. Here and throughout,
$\log_j$ stands for the $j$-fold iterated logarithm. Later on, these
estimates have been improved: Corr\'adi and K\'atai [2] obtained $$
P(x) = \Omega_+\(x^{1/2}\,\exp\(C_1(\log_2 x)^{1/4}\,(\log_3
x)^{-3/4}\)\)\,, \eqno(2.5)$$ Hafner [8] established $$
P(x)=\Omega_-\(x^{1/2}\,(\log x)^{1/4}(\log_2
x)^{(\log2)/4}\,\exp(-C_2(\log_3 x)^{1/2})\)\,, \eqno(2.6)$$ and
Soundararajan [22] proved that $$ P(x)=\Omega\(x^{1/2}\,(\log
x)^{1/4}(\log_2 x)^{3(2^{1/3}-1)/4}\,(\log_3 x)^{-5/8}\)\,.
\eqno(2.7)$$ The bounds (2.5) - (2.7) depend on the special
multiplicative structure of the arithmetic function $r(n)$, and on
the analytic properties of its generating Dirichlet series (Epstein
zeta-function). \bsk\bsk

{\bf 3.~Results on the spectral counting function of Heisenberg
manifolds. } Returning to rational Heisenberg manifolds
$M=(H_\l/\Gamma_{\b r}, g_\l)$ as described in section 1, we give an
account of what is known about the error term in (1.1), i.e.,
$$ R(t) = N(t)-
{\vol(M)\over(4\pi)^{\l+1/2}\,\Gamma(\l+{3\over2})}\,t^{\l+1/2}=
N(t) - {\rr\over2^{2\l+1/2}\pi^{\l}\,
\Gamma(\l+{3\over2})}\,t^{\l+1/2}\,. \eqno(3.1) $$ For $\l=1$,
Petridis and Toth [21] proved that $R(t)\ll t^{5/6}\log t$. This
estimate was sharpened and generalized to arbitrary $\l\ge1$ by
Khosravi and Petridis [16] who established $R(t)\ll t^{\l-7/41}$. In
a recent paper, Zhai [24] applied Huxley's "discrete
Hardy-Littlewood method" [11], [12] to derive, for any $\l\ge1$,
$$ R(t) \ll t^{\l-77/416}(\log t)^{26957/8320}\,. \eqno(3.2)$$ Actually,
it is just the special "rational" choice of the metric $g_\l$ which
makes the error term (possibly) large. As Khosravi and Petridis [16]
showed, for "almost all" metrics $g$, the much sharper bound
$$ R_g(t)\ll_g\ t^{\l-1/4}\,\log t \eqno(3.3) $$ holds true.
Returning to the rational case (1.4), a result of Khosravi [15]
and Khosravi \& Toth [17] tells us that $$ \Int_0^T (R(t))^2 \d t
= C_\l\,T^{2\l+1/2} + \O{T^{2\l+1/4+\ep}} \eqno(3.4) $$ where
$C_\l>0$ is an explicit constant. A recent paper of Zhai [24] is
concerned with estimates and asymptotics for higher power moments
of $R(t)$. In fact, (3.3) and (3.4) may suggest the conjecture
that
$$ R(t) \ll t^{\l-1/4+\ep} \eqno(3.5) $$ for every $\ep>0$.
The results described so far show a lot of analogy to the Gaussian
circle problem discussed in section 2. \ssk

In the present work, it is our objective to estimate $R(t)$ from
below, in order to arrive again at a statement saying that \
"$R(t) \ll t^{\l-1/4}$ {\it in mean-square, with an unbounded
sequence of exceptionally large values $t$}". In fact, we are able
to find for each $\l\ge1$ an explicit function $\om_\l(t)$ tending
to $\infty$, such that $$ R(t)=\Omega(t^{\l-1/4}\om_\l(t))\,.  $$
\bsk

\bsk {\bf Theorem. } { \it For any fixed positive integer $\l$, let
$(H_\l/\Gamma_{\b r},g_\ell)$ be a rational $(2\l+1)$-dimensional
Heisenberg manifold with metric $g_\l$, as described above. Then the
error term $R(t)$ for the associated spectral counting function,
defined in $(3.1)$, satisfies
$$ \limsup_{t\to\infty}{R(t)\over t^{\l-1/4}\,\om_\l(t)}> 0\,, $$ where
$$ \om_\l(t):=\cases{(\log t)^{1/4}& for $\l$ even,\cr
(\log_2 t\,\log_3 t)^{1/4}& for $\l$ odd.\cr} $$ } \bsk

{\bf Remarks. } 1. The case of even $\l$ has been treated in the
first part of this work [20]. After approximating $R(t)$ by a
suitable trigonometric sum, the Dirichlet approximation theorem
was applied to give all its terms the positive sign. In the
present article we shall deal with the case of odd $\l$, employing
a quantitative version of Kronecker's theorem instead.
Technically, this will be stated in terms of uniform distribution
theory - see Lemma 4 below. \ssk 2. Our results obviously are
comparable to the bounds (2.3) and (2.4) for the circle problem.
It seems very difficult to obtain improvements as sharp as (2.5) -
(2.7), since the coefficients $\theta_\l(n)$ (defined in (5.11)
below) fail to share the useful properties of $r(n)$. \ssk 3. For
the circle problem, the two different types of arguments
(Dirichlet's theorem vs. Kronecker's) were used to establish
$\Omega_-$- and $\Omega_+$-results. For Heisenberg manifolds they
are needed to deal with $\l$ of arbitrary parity, yielding
$\Omega_+$-bounds in both cases.

\bsk\bsk

%\bye

\bsk\bsk

{\bf 4.~Some Lemmas.}\bsk

{\bf Lemma 1. } (Vaaler's approximation of fractional parts by
trigonometric polynomials.) { \it For arbitrary $w\in\R$ and
$H\in\Z^+$, let $\psi(w):= w-[w]-\h$, $$ \Sigma_H(w):=\sum_{h=1}^H
{\al_{h,H}}\,\sin(2\pi hw)\,,\qquad \Sigma_H^*(w):=\sum_{h=1}^H
{\be_{h,H}}\,\cos(2\pi hw)\ + {1\over2H+2}\,, $$ where, for
$h=1,\dots,H$,
$$ \al_{h,H}:={1\over\pi h}\,\rho\({h\over H+1}\)\,,\qquad \be_{h,H}:={1\over H+1}\(1-{h\over
H+1}\)\,, $$ and $$ \rho(\xi)= \pi\xi(1-\xi)\cot(\pi\xi)+\xi
\qquad\qquad (0<\xi<1)\,. $$ Then the following inequality holds
true:
$$ \abs{\psi(w)+\Sigma_H(w)} \le \Sigma_H^*(w)\,. $$ } \bsk

{\bf Proof. } This is one of the main results in Vaaler [23]. A very
well readable exposition can also be found in the monograph by
Graham and Kolesnik [7].

 \bsk\msk

{\bf Lemma 2. } { \it Let $F\in C^4[A, B]$,  $G\in C^2[A, B]$, and
suppose that, for positive parameters $X, Y, Z$, we have $1\ll
B-A\ll X$ and
$$ F^{(j)}\ll X^{2-j} Y^{-1} \for j=2,3,4, \
\abs{F''}\ge c_0 Y^{-1}\,,\quad G^{(j)}\ll X^{-j} Z \for j=0,1,2,
$$ throughout the interval $[A,B]$, with some constant $c_0>0$.
Let $\J'$ denote the image of $]A,B]$ under $F'$, and $F^*$ the
inverse function of $F'$. Then, with $e(w)=e^{2\pi iw}$ as usual,
$$ \eq{\sum_{A<m\le B} G(m)\,e(F(m)) =&\ e\({{\rm sgn}(F'')\over8}\) \sum_{k\in\J'}{G(F^*(k))
\over\sqrt{\abs{F''(F^*(k))}}}\,e\(F(F^*(k))-kF^*(k)\)  + \cr & +
\O{Z\(\sqrt{Y}+\log(2+\,{\rm length}(\J'))\)}\,. \cr }
$$ }\bsk

{\bf Proof. } Transformation formulas of this kind are quite common,
though often with worse error terms. This very sharp version can be
found as f.~(8.47) in the recent monograph [14] of H.~Iwaniec and
E.~Kowalski.

\bsk\msk

{\bf Lemma 3. } {\it For a real parameter $T\ge1$, let $\Fe_T$
denote the Fej\'er kernel $$ \Fe_T(v) = T\({\sin(\pi Tv)\over\pi
Tv }\)^2\,. $$ Then for arbitrary real numbers $Q>0$ and $\delta$,
it follows that $$ \Int_{-1}^1 \Fe_T(v) e(Qv+\delta)\d v =
\max\(1-{Q\over T},0\)e(\delta) + \O{1\over Q}\,, $$ where the
$O$-constant is independent of $T$ and $\delta$.}\bsk

{\bf Proof. } This useful result is due to Hafner [8]. It follows
from the classic Fourier transform formula $$ \Int_\Ri \Fe_T(v)
e(Qv) \d v = \Int_\Ri \({\sin(\pi v)\over\pi v }\)^2\ e\({Q\over
T}v\)\d v = \max\(1-{Q\over T},0\)\,. $$ Since $\Fe_T(\pm1)\ll
T^{-1}$ and $\Fe_T'(v)\ll v^{-2}$ for $|v|\ge1$, uniformly in
$T\ge1$, integration by parts readily shows that the intervals
$]-\infty,-1]$ and $[1,\infty[$ contribute only $\O{Q^{-1}}$.
\bsk\msk

{\bf Lemma 4. } { \it For an arbitrary integer $s\ge2$, let $\b a=
(a_1,\dots,a_s)\in\R^s$ so that $1,a_1,\dots,a_s$ are linearly
independent over $\Z$. Suppose further that there exists a function
$\phi:\R^+\to\R^+$ such that $\phi(t)/t$ increases monotonically and
$$ \norm{\b h\cdot\b a}\ge{1\over\phi(|\b h |_\infty)}$$ for all $\b
h\in\Z^s\setminus\{\b o\}$, where $\norm{\cdot}$ denotes the distance from
the nearest integer. Then for any positive integers $N_0$ and $N$,
the discrepancy modulo 1 $D_{N_0,N}(n\b a)$ of the sequence $\(n\b a
\)_{n=N_0+1}^{N_0+N}$ satisfies $$ D_{N_0,N}(n\b a) \le
c^s\,s!\,{\log N \log\phi^{-1}(N)\over\phi^{-1}(N)}\,,$$ where $c$
is an absolute constant, $\phi^{-1}$ denotes the inverse function of
$\phi$, and $N$ is supposed to be so large that $\phi^{-1}(N)\ge
e$.}\bsk\msk

{\bf Proof. } This is essentially Theorem 1.80 in the monograph of
Drmota and Tichy [4], p.~70, with the dependance on the dimension
worked out explicitly.

\bsk\bsk {\bf 5.~Proof of the Theorem. }\quad As already stated, the
case of even $\l$ has been treated in part I of this work [20].
Therefore, we may suppose throughout that $\l$ is odd. We start from
Lemma 3.1 in Zhai [24] which approximates the error term involved by
a fractional part sum. Let $U$ be a large real parameter,
$u\in[U-1,U+1]$, and put
$$ E(u):= {2^{\l-2}(\l-1)!\over\rr}\,R(2\pi u^2)\,.
\eqno(5.1) $$ Then according to Zhai\fn{5}{In fact, Zhai in his
notation tacitly assumes that $r_1=\dots=r_\l=1$, which means no
actual loss of generality. We have supplemented the factor $\rr$ in
(5.1).} [24], Lemma 3.1, for arbitrary\fn{6}{At this stage, we write
up the argument for general $\l$, in order to point out the
importance of the condition that $\l$ is odd later on.} $\l\ge1$,
$$ \eq{E(u) &= E^*(u) + \O{u^{2\l-1}}\,, \cr E^*(u) &:= -
\sum_{1\le m\le u}
m(u^2-m^2)^{\l-1}\psi\({u^2\over2m}-{m\over2}-{\l\over2}\) \,.\cr}
\eqno(5.2)
$$ We apply Lemma 1 in the form $-\psi\ge\Sigma_H - \Sigma_H^*$,
choosing $H=[U]$. Thus we get $$ \eq{&E^*(u)\ge \ - U  + \DS
m(u^2-m^2)^{\l-1}\times\cr &\times(-1)^{h\l}\(\alh\sin\(2\pi
h\({u^2\over2m}-{m\over2}\)\)- \beh\cos\(2\pi
h\({u^2\over2m}-{m\over2}\)\)\) \,.\cr} \eqno(5.3)$$ We split up the
range $1\le m\le u$ into dyadic subintervals $\Mj=]M_{j+1},M_j]$,
$M_j=u\,2^{-j}$ for $j=0,\dots,J$, where $J$ is minimal such that
$(U-1)2^{-J-1}<1$. We thus have to deal with exponential sums $$
\E_j(h,u):=\sum_{m\in\Mj}m(u^2-m^2)^{\l-1}\,
e\(-h\({u^2\over2m}-{m\over2}\)\)\,. $$ We transform them by Lemma
2, with $ G(\xi) = \xi(u^2-\xi^2)^{\l-1}$, $F(\xi)=
-h\({u^2\over2\xi}-{\xi\over2}\)$. On each interval $\Mj$ the
conditions of Lemma 2 are fulfilled with the parameters $X= M_j$,
$Y={M_j^3\over hu^2}$, $Z=M_j\,u^{2\ell-2}$. By straightforward
computations, as in [20], we obtain $$ \eq{\E_j(h,u) = & \
h^{3/4}\,u^{2\l-1/2} \sum_{k\in F'(\Mj)}
{(2k-2h)^{\l-1}\over(2k-h)^{\l+1/4}}\,e\(-u\sqrt{h}\sqrt{2k-h}-{\textstyle{
1\over8}}\)\cr &+ \O{u^{2\l-3}{M_j^{5/2}\over h^{1/2}}+u^{2\l-1}\log
u}\,.\cr}  \eqno (5.4)$$ It is plain to see that the overall
contribution of the error terms to (5.3), summed over $j$ and $h$,
is $\ll u^{2\l-1/2}$: See [20], f.~(5.5). Summing up the main terms
in (5.4), the total range of $k$ becomes $h=F'(M_0)< k \le
F'(M_{J+1})=\h h+2^{2J+1} h=:K_{h,U}$, hence
$$ \eq{&\DS m(u^2-m^2)^{\l-1} (-1)^{h\l}\gah
\,e\(-h\({u^2\over2m}-{m\over2}\)\) = \cr &= u^{2\l-1/2} \Sh
(-1)^{h\l}\gah\,h^{3/4}\, \sum_{h<k\le K_{h,U}}
{(2k-2h)^{\l-1}\over(2k-h)^{\l+1/4}}\,e\(-u\sqrt{h}\sqrt{2k-h}-{\textstyle{
1\over8}}\) + \cr & + \O{u^{2\l-1/2}}\,,\cr} \eqno(5.5) $$ where
$\gah$ stands for either $\alh$ or $\beh$. Using the real and
imaginary part of this result in (5.3), we arrive at
$$ E^*(u) \ge u^{2\l-1/2}\,S(u,U) - c_1 u^{2\l-1/2}\,, \eqno(5.6) $$
where $$ \eq{&S(u,U) :=\ \sum_{(h,k)\in\D(U)} h^{3/4} \,
{(2k-2h)^{\l-1}\over(2k-h)^{\l+1/4}}\times\cr&\times(-1)^{h\l}\(\alh\,
\sin\(2\pi u\sqrt{h}\sqrt{2k-h}+{\textstyle{\pi\over4}}\)-\beh\,
\cos\(2\pi u\sqrt{h}\sqrt{2k-h}+{\textstyle{\pi\over4}}\)\)\,,\cr}
$$ $$ \D(U) := \{(h,k)\in\Z^2:\ 1\le h\le U\,,\ h<k\le
K_{h,U}\ \}\,, $$ and $c_1$ is an appropriate positive constant.
The next step is to eliminate the majority of the terms of the
last double sum. To this end, let $T$ be another large parameter,
with the constraint that $U\ge T^2$. Using Lemma 3, we multiply
$S(u,U)$ by the Fej\'er kernel $\Fe_T(u-U)$ and integrate over
$U-1\le u\le U+1$. Thus
$$ \eq{& I(T,U) := \Int_{U-1}^{U+1} S(u,U)\Fe_T(u-U)\d u =
\Int_{-1}^1 S(U+v,U)\Fe_T(v)\d v = \cr & = \sum_{(h,k)\in\D(U),\
h(2k-h)\le T^2} h^{3/4} \,
{(2k-2h)^{\l-1}\over(2k-h)^{\l+1/4}}\(1-{\sqrt{h(2k-h)}\over T
}\)(-1)^{h\l}\times\cr&\times\(\alh\, \sin\(2\pi
U\sqrt{h(2k-h)}+{\textstyle{\pi\over4}}\)-\beh\, \cos\(2\pi
U\sqrt{h(2k-h)}+{\textstyle{\pi\over4}}\)\)\cr & +
\O{\sum_{(h,k)\in\D(U)} h^{-3/4}
\,{(2k-2h)^{\l-1}\over(2k-h)^{\l+3/4}}}\,.\cr } \eqno(5.7) $$ The
$O$-term here is in fact $O(1)$: see [20], f.~(5.9). Now recall that
$(h,k)\in\D(U)$ explicitly means that $$ 1\le h\le U\,,\quad h<k\le
\h h+2^{2J+1} h \asymp U^2 h\,, $$ while $h(2k-h)\le T^2\,,\ k>h$
implies that $$ h\le T\le\sqrt{U}\quad\and\quad k<2k-h\le {T^2\over
h} \le T^2\le U\,.$$ Hence the summation condition on the right hand
side of (5.7) can be simplified to $$ 0<h(2k-h)\le T^2\,,\quad
k>h\,. \eqno(5.8) $$ For any $h, k$ satisfying (5.8), write
$h(2k-h)=r^2 q$, with $r$ an integer and $q$ square-free. Now
suppose we can choose $U$ so that $$ \norm{U\sqrt{q}-{1\over2}}\le
{\ep_0\over T} \eqno(5.9) $$ for all square-free $q\in]1,T^2]$. Here
$\norm{\cdot}$ denotes the distance from the nearest integer and $\ep_0>0$
is a suitably small constant. Then, for $r^2 q\le T^2$, $q>1$, $$
\sin\(2\pi U\,r\sqrt{q}+{\textstyle{\pi\over4}}\)=
{(-1)^{r}\over\sqrt{2}}+\O{\ep_0}\,. $$ By the definitions in Lemma
1, for $h\le T\le\sqrt{U}$, it follows that $\alh={1\over\pi
h}+o(1)$, $\beh\ll{1\over U}$. Hence, provided that $(5.9)$ is true,
for $h(2k-h)$ not a perfect square, $$ \eq{&\alh\, \sin\(2\pi
U\sqrt{h(2k-h)}+{\textstyle{\pi\over4}}\)-\beh\, \cos\(2\pi
U\sqrt{h(2k-h)}+{\textstyle{\pi\over4}}\)\cr & =\ {1\over h
}\({(-1)^{h}\over\pi\sqrt{2}}+\O{\ep_0}+o(1)\)\,,\cr} \eqno(5.10)
$$ where $o(1)$ refers throughout to $U\to\infty$. (Note that $h$,
$h(2k-h)=r^2q$, and $r$ are of the same parity.) We can use (5.10)
in (5.7), appealing at last to the condition that $\l$ is odd: Then
$(-1)^{h\l}(-1)^{h}=1$, i.e., the alternating factors cancel out. To
simplify notation, we put $$ \theta_{\l}(n) := \sum_{h(2k-h) = n,\
k>h} {h^{1/2}\over(2k-h)^{1/2}}\,\(1-{h\over
2k-h}\)^{\ell-1}\,.\eqno(5.11)
$$ With that and (5.10), eq.~(5.7) readily yields\fn{7}{Properly
speaking, we are committing here an error concerning the $n$ which
are perfect squares. But it is obvious that their contribution to
the right-hand side of (5.12) is $\ll\sum m^{-3/2+\ep}\ll1$.}
$$ \eq{I(T,U) &\ge c_2\sum_{1\le n\le T^2} {\theta_{\l}(n)\over
n^{3/4}} \(1-{\sqrt{n}\over T}\)\ - c_3\cr &\ge c_4 \sum_{1\le n\le
T^2/2} {\theta_{\l}(n)\over n^{3/4}} \ - c_3\,, \cr}\eqno(5.12) $$
with certain positive constants $c_2, c_3, c_4 $. Further, for $1\le
n\le T^2$, $$ \theta_\l(n)\gg \sum_{h(2k-h)=n\atop k>2h}
{h^{1/2}\over(2k-h)^{1/2}} \gg \sum_{hm=n,\ m>3h\atop h\equiv m
\mod2}{\sqrt{h}\over\sqrt{m}}\gg \sum_{hm=n,\ 3h<m\le4h\atop h\equiv
m \mod2}1\,.   $$ Therefore, by (5.12), $$ \eq{I(T,U) + c_3 &\gg
\sum_{T^2/4\le n\le T^2/2 \atop n\odd} n^{-3/4} \sum_{hm=n\atop
3h<m\le4h } 1 \cr &\gg\ T^{-3/2} \sum_{T^2/4\le hm\le T^2/2\atop
3h<m\le4h,\ h,m \odd}1 \ \gg\ T^{1/2}\,.\cr} \eqno(5.13)$$ It
remains to ensure the validity of (5.9) and, at the same time,
establish a lower bound for $T$ in terms of $U$. To this end we will
employ Lemma 4. Let $q_1,\dots,q_s$ be all square-free integers in
$]1,T^2]$, then obviously $s\asymp T^2$. Let further
$p_1<p_2<\dots<p_P$ be all primes \hbox{$\le T^2$,} with
$P=\pi(T^2)\asymp T^2/\log T$. Consider the field extension
$\F=\Q(\sqrt{p_1},\dots,\sqrt{p_P})$ and the corresponding Galois
group $G={\rm Gal}(\F/\Q)$. Since $\chi(\sqrt{p_j})=\pm\sqrt{p_j}$
for $j=1,\dots,P$ and every $\chi\in G$, obviously $|G|=2^P$.\ssk
According to Besicovitch's theorem [1], the numbers
$1,\sqrt{q_1},\dots,\sqrt{q_s}$ are linearly independent over $\Z$.
For any $\b h=(h_1,\dots,h_s)\in\Z^s\setminus\{\b o\}$, let $-h_0$
denote the integer nearest to $\b
h\cdot(\sqrt{q_1},\dots,\sqrt{q_s})$, then
$$ \abs{\prod_{\chi\in
G}\chi(h_0+h_1\sqrt{q_1}+\dots+h_s\sqrt{q_s})}\ge\ 1\,, $$ since the
left-hand side is the absolute value of the norm of a nonzero
algebraic integer. For any $\chi\in G$, $$
\abs{\chi(h_0+h_1\sqrt{q_1}+\dots+h_s\sqrt{q_s})}\le3s|\b h|_\infty
T\,. $$ Hence, $$
\abs{h_0+h_1\sqrt{q_1}+\dots+h_s\sqrt{q_s}}\ge\(3s|\b h|_\infty
T\)^{-2^P}\,. $$ Therefore, the requirements of Lemma 4 are
fulfilled with $\phi(t)=(3s\,t\,T)^{2^P}$. Thus we may conclude
that, for $N$ sufficiently large, the discrepancy modulo 1 $D_N$ of
the sequence $\(n\sqrt{q_1},\dots,n\sqrt{q_s}\)_{n=N+1}^{2N}$
satisfies $$ D_N \le c^s\,s!\,{\log N
\log\phi^{-1}(N)\over\phi^{-1}(N)} \ll c^s s!\,N^{-2^{-P}}(\log N)^2
\,,$$ since $\phi^{-1}(N)=N^{2^{-P}}(3sT)^{-1}$. Now, with (5.9) in
the back of mind, let $\c$ denote the $s$-dimensional cube
$[\he-{\ep_0\over T},\he+{\ep_0\over T}]^s$. Writing
$\langle\cdot\rangle$ for the fractional part, we may thus conclude
that $$ \eq{&\#\{n\in\Z:\ N<n\le2N,\ \(\langle
n\sqrt{q_1}\rangle,\dots,\langle n\sqrt{q_s}\rangle\)\in\,\c\ \}\cr
&\ge N\vol(\c)-ND_N\ge N\((2\ep_0)^s T^{-s}-c_5^s
s!\,N^{-2^{-P}}(\log N)^2\)\cr &\ge N\((c_6T)^{-c_7T^2} -
(c_8T)^{c_9T^2} N^{-2^{-P}}(\log N)^2\)\,,  \cr} $$ using Stirling's
formula and $s\asymp T^2$. Here $c_j$ are throughout absolute
positive constants. A short calculation shows that the large bracket
is certainly $>0$ if we choose $$ N=N^*(T):=[\exp(\exp(c_{10}T^2/
\log T))]\,, \eqno(5.14)  $$ with $c_{10}$ sufficiently large.
Therefore, for arbitrary $T$ there exists at least one integer
$U=n\in\ ]N^*(T), 2N^*(T)]$ for which $(\langle
U\sqrt{q_1}\rangle,\dots,\langle U\sqrt{q_s}\rangle)\in \c$, which
means that (5.9) is fulfilled. Further, (5.14) implies that $$ T
\asymp (\log_2 U\,\log_3 U)^{1/2}\,. \eqno(5.15) $$ Thus (5.13)
yields $$ I(T,U) \gg (\log_2 U\,\log_3 U)^{1/4}\,. \eqno(5.16) $$ On
the other hand, it follows from the definition of $I(T,U)$ that $$
I(T,U) \le \(\sup_{U-1\le u\le U+1 }S(u,U)\)\Int_{-1}^1\Fe_T(v)\d
v\le \sup_{U-1\le u\le U+1 }S(u,U)\,. $$ This implies that there
exists a value $u^* \in [U-1,U+1]$ for which $$ S(u^*,U)\gg (\log_2
u^*\,\log_3 u^*)^{1/4}\,. \eqno(5.17) $$ It remains to recall that
if $T$ runs through an unbounded sequence of positive reals, by
construction\fn{8}{Recall the ultimate order in the choice of the
parameters: For $T$ an independent large real variable, $N^*(T)$ is
defined by (5.14), then $U$ is picked from $]N^*(T),2N^*(T)]$,
finally $u^*$ from $[U-1,U+1]$. The condition $U\ge T^2$ needed
earlier in the argument is amply satisfied in view of (5.15).} so do
$U$ and $u^*$. Therefore, (5.1), (5.2), (5.6) and (5.17) together
complete the proof of our theorem.

\bsk  \bsk

{\klein  \parindent=0pt \def\smc{}

\cen{\bf References}  \bsk %\hsize=17true cm \vsize=24.5true cm

[1] A.S.~Besicovitch, On the linear independence of fractional
powers of integers, J.~London Math. Soc. {\bf15}, 3-6 (1940).\ssk

[2] {K.~Corr\'adi and I.~K\'atai,} A comment on K. S. Gangadharan's
paper "Two classical lattice point problems" (Hungarian), Magyar
Tud. Akad. mat. fiz. Oszt. K\"ozl. {\bf 17}, 89-97 (1967).  \ssk

[3] H.~Cram\'er, \"{U}ber zwei S\"{a}tze von Herrn G.H.~Hardy, Math.~Z.
{\bf15}, 201-210 (1922). \ssk

[4] M.~Drmota and R.F.~Tichy, Sequences, discrepancies, and
applications, Lecture Notes in Math. 1651, Springer, Berlin, 1997.
\ssk

[5] {K.S.~Gangadharan,} Two classical lattice point problems,
Proc.~Cambridge Phil.~Soc. {\bf 57}, 699--721 (1961).  \ssk

[6] C.S.~Gordon and E.N.~Wilson, The spectrum of the Laplacian on
Riemannian Heisenberg manifolds, Michigan Math. J. {\bf33}, 253-271
(1986). \ssk

[7] S.W.~Graham and G.~Kolesnik, Van der Corput's method of
exponential sums, Cambridge 1991. \ssk

[8] J.L. Hafner, New omega results for two classical lattice point
problems, Invent. Math. {\bf63}, 181-186 (1981).\ssk

[9] {G.H.~Hardy,} On the expression of a number as the sum of two
squares, Quart. J.~Math. {\bf 46}, 263--283 (1915). \ssk

[10] L.~H\"{o}rmander, The spectral function of an elliptic operator,
Acta Math. {\bf121}, 193-218 (1968).\ssk

[11] {\smc M.N.~Huxley}, {Area, lattice points, and exponential
sums,} LMS Monographs, New Ser. {\bf 13}, Oxford 1996.  \ssk

[12] M.N.~Huxley, Exponential sums and lattice points III.
Proc.~London Math.~Soc. (3) {\bf87}, 591-609 (2003). \ssk

[13] {\smc A.~Ivi\'c, E.~Kr\"{a}tzel, M.~K\"{u}hleitner, \and W.G.~Nowak,}
Lattice points in large regions and related arithmetic functions:
Recent developments in a very classic topic. Proceedings Conf.~on
Elementary and Analytic Number Theory ELAZ'04, held in Mainz, May
24-28, W.~Schwarz and J.~Steuding eds., Franz Steiner Verlag 2006,
pp. 89-128.\ssk

[14] {\smc H.~Iwaniec, E.~Kowalski,} Analytic Number Theory, AMS
Coll.Publ.~53. Providence, R.I., 2004. \ssk

[15] M.~Khosravi, Spectral statistics for Heisenberg manifolds,
Ph.D.~thesis, McGill U. 2005. \ssk

[16] M.~Khosravi and Y.N.~Petridis, The remainder in Weyl's law for
$n$-dimensional Heisenberg manifolds, Proc. AMS {\bf133}/12,
3561-3571 (2005).\ssk

[17] M.~Khosravi and J.A.~Toth, Cramer's formula for Heisenberg
manifolds, Ann.~de l'institut Fourier {\bf55}, 2489-2520 (2005).\ssk

[18] {\smc E.~Kr\"atzel,} Lattice points. Berlin 1988. \ssk

[19] {\smc E.~Kr\"atzel,} Analytische Funktionen in der
Zahlentheorie. Stuttgart-Leipzig-Wiesbaden 2000.  \ssk

[20] W.G.~Nowak, A lower bound for the error term in Weyl's law for
certain Heisenberg manifolds, submitted for publication. Available
online at: {\tt
http://arxiv.org/PS\_cache/arxiv/pdf/0809/0809.3924v1.pdf} \ssk

[21] Y.N.~Petridis and J.A.~Toth, The remainder in Weyl's law for
Heisenberg manifolds, J.~Diff.~Geom. {\bf60}, 455-483 (2002).\ssk

[22] {\smc K.~Soundararajan,} Omega results for the divisor and
circle problems. Int.~Math.~Res.~Not. {\bf36}, 1987-1998 (2003).
\ssk

[23] J.D.~Vaaler, Some extremal problems in Fourier analysis,
Bull.~Amer.~Math.~Soc. {\bf12}, 183-216 (1985).\ssk

[24] W.~Zhai, On the error term in Weyl's law for the Heisenberg
manifolds, Acta Arithm.~{\bf134}, 219-257 (2008).

\vbox{\vskip 0.7true cm}

\parindent=1.5true cm

\vbox{Institute of Mathematics

Department of Integrative Biology

Universit\"at f\"ur Bodenkultur Wien

Gregor Mendel-Stra\ss e 33

1180 Wien, \"Osterreich \ssk

E-mail: {\tt nowak@boku.ac.at} \ssk

}}

\bye